\documentclass[a4paper]{article}

\usepackage{
amsmath,
amsthm,
amscd,
amssymb,
}
\usepackage{comment}
\usepackage{tikz}
\usetikzlibrary{positioning,arrows,calc}
\usepackage{xspace}

\usepackage{graphicx}

\usepackage{url}

\theoremstyle{plain}

\newtheorem{definition}{Definition}

\newtheorem{theorem}{Theorem}

\newtheorem{question}{Question}

\newtheoremstyle{derp}
{3pt}
{3pt}
{}
{}
{\upshape}
{:}
{.5em}
{}
\theoremstyle{derp}

\newcommand{\Z}{\mathbb{Z}}

\newcommand{\B}{\mathcal{L}}

\title{Minimal subshifts with a language pivot property}

\author{
Ville Salo \\
vosalo@utu.fi
}

\begin{document}
\maketitle

\begin{abstract}
We construct a binary minimal subshift whose words of length $n$ form a connected subset of the Hamming graph for each $n$.
\end{abstract}

\section{Background}

We paraphrase a question M. Hochman asked in the conference Current Trends in Dynamical Systems and the Mathematical Legacy of Rufus Bowen, organized in Vancouver in 2017. This question (and many more) can be found on the website of the conference \cite{Trends}.

\begin{question}
\label{q}
Does there exist a minimal subshift $X$ such that for some $k$ and all $n$, the set of words of length $n$ in $X$ is connected in the $k$th power of the Hamming graph? Is such $X$ necessarily topologically mixing?
\end{question}

We give two solutions: For $k = 2$, we prove Sturmian subshifts are an example. For $k = 1$, we construct such a minimal subshift from scratch. We do not solve the second question for $k = 1$. 

The example for $k = 1$ was written right after the conference. I have fixed some typos later. The Sturmian example arose in discussions with Nishant Chandgotia, who gave another proof of Theorem~\ref{thm:Sturmian} in \cite{chandgotia2018pivot}.

\section{Definitions}

An \emph{alphabet} is a finite set $\Sigma$ of \emph{symbols}. A \emph{(formal) language} (over $\Sigma$) is a subset of $\Sigma^*$, where $\Sigma^*$ is the set of all (possibly empty) words. We write $\epsilon$ for the empty word. Words over $\Sigma$ form a free monoid with basis $\Sigma$ and identity element $\epsilon$. The monoid operation is concatenation, and is written as $u \cdot v$ or simply $uv$.

We say a set of words $W \subset \Sigma^n$ is \emph{connected} if the graph with nodes $W$ and edges $(u, v)$ for all $u, v \in W$ having Hamming distance $H(u, v) = 1$ is connected.

If $X \subset \Sigma^\Z$ is a subshift (closed shift-invariant set), we write $\B(X) \subset \Sigma^*$ for its \emph{language} (set of words that appear in its points) and $\B_n(X) = \B(X) \cap \Sigma^n$. We say $X$ is \emph{language-connected} if for every $n$ the language $\B_n(X)$ is connected. A subshift is \emph{minimal} if it has no nontrivial subshifts, equivalently every word in the language appears as a subword of every long enough word in the language. A subshift $X$ is \emph{topologically mixing} (sometimes called strongly mixing) if for any two words $u, v$ in the language of $X$, there exists $n$ such that
\[ m \geq n \implies \exists w: |w| = m \wedge uwv \in \B(X). \]

If $L \subset \Sigma^*$ is a language which is \emph{extendable} in the sense that
\[ v \in L \implies \exists u, w: |u|, |w| > 0 \wedge uvw \in L, \]
then we write $\B^{-1}(L)$ for the subshift $X \subset \Sigma^\Z$ whose language is the closure of $L$ under taking subwords, that is,
\[ \B(X) = \{v \in \Sigma^* \;|\; \exists u, w: uvw \in L\} \]

We use the usual semiring structure for languages, in particular $L \cdot K = \{u \cdot v \;|\; u \in L, v \in K\}$. We also use the Kleene star operation for $A \subset \Sigma^*$:
\[ A^* = \{w_1 w_2 \cdots w_k \;|\; k \geq 0, \forall i \in \{1,\ldots,k\}: w_i \in A\}. \]

\section{Sturmian subshifts}

It turns out that for Sturmian subshifts (which are minimal and not strongly mixing), the subshift has a connected language when two changes are allowed at once, and we can even make the two changes next to each other. Thus, Sturmian subshifts solve Question~\ref{q}.

\begin{definition}
The \emph{$k$-change graph} of a language $L \subset \Sigma^n$ has nodes $L$ and edge $(u, v)$ whenever 
\[ u = wu'w', v = wv'w' \]
where $u, u' \in \Sigma^{k'}$ for some $k' \leq k$ and $w,w' \in \Sigma^*$.
\end{definition}

\begin{theorem}
\label{thm:Sturmian}
Let $X$ be a Sturmian subshift. Then for every $n \geq 2$, the $2$-change graph of $\B_n(X)$ is connected.
\end{theorem}

\begin{proof}
Configurations in a Sturmian subshift can be characterized as mechanical words, i.e. for a Sturmian subshift there exists an irrational slope $\alpha > 0$ such that every configuration is obtained by drawing a line of slope $\alpha$ on the standard planar embedding of the Cayley graph of $\Z^2$ (i.e. a regular square grid), and writing $0$ whenever the line crosses a horizontal line and $1$ whenever it crosses a vertical line (and the line does not contain an element of $\Z^2$).

In particular, words of length $n$ in a Sturmian subshift are obtained by starting a line (that does not contain elements of $\Z^2$) from some point on $A = ((0,1) \times \{0\}) \cup (\{0\} \times (0,1))$ and following it until it crosses $n$ horizontal or vertical lines.

Thus if $u, v \in \B_n(X)$ and $X$ is Sturmian, then both $u$ and $v$ have corresponding finite line segments that start on some points $x_a, x_b \in A$ and end on some horizontal or vertical line of the grid. We think of the line as crossing the horizontal or vertical line of $A$ that it begins on by an infinitesimal amount (depending on the first symbol of $u$), and similarly the line on which it ends is crossed infinitesimally, determined by the last symbol.

Now slide $x_a$ towards $x_b$ along $A$ and have the path cross exactly $n$ lines at all times. The word coded by crossings varies continuously (i.e. does not change at all) as long as the line does not hit an element of $\Z^2$. There are finitely many points where the line does hit an element of $\Z^2$. By irrationality of $\alpha$, it cannot hit two at once, so we see that when a central point of the line hits an element of $\Z^2$, we will only change $01$ to $10$ or vice versa in the word. When we move from $(0,1) \times \{0\}$ to $\{0\} \times (0,1)$ on $A$, we only flip the initial symbol of the word. When the endpoint of the line crosses an integer point, either only the last symbol changes or the last two symbols change.
\end{proof}

I do not know whether it is possible to somehow turn a Sturmian subshift into a language-connected minimal subshift. No simple recoding trick seems to work.

\section{Language-connected example}

\begin{theorem}
There exists a language-connected minimal subshift.
\end{theorem}

\begin{proof}
Let $\Sigma = \{0,1\}$. Let $w_{0,0} = 01$ and $w_{0,1} = 11$, $W_0 = \{w_{0,0}, w_{0,1}\}$, $\ell_0 = 2$, $n_0 = 2$.

Assuming $W_i$ has been defined, define
\[ L_{i,k,g} = (W_i \{\epsilon, 0, \ldots, 0^g\})^{k-1} W_i. \]
Note that $L_{i,k,0} = W_i^k$. If $w$ is a word over this language and we have fixed such a decomposition into the defining form, then we call the maximal words $0^a$ separating words in $W_i$ \emph{gaps}, and refer to the words taken from $W_i$ as the \emph{$W_i$-words} of $w$. There may be multiple decompositions of a word of $L_{i,k,g}$, but we usually work with, and consistently modify, a fixed such decomposition. This should not cause confusion.

For each $i$, we build by induction a set of $n_i$ words $W_i = \{w_{i,j} \;|\; j \in \{0, \ldots, n_i-1\}\} \subset \Sigma^{\ell_i}$ with the following properties:
\begin{itemize}
\item $H(w_{i,j}, w_{i,j+1}) = 1$ for all $j \in \{0,\ldots,n_i-2\}$.
\item $H(w_{i,0} 0, 0 w_{i,n_i-1}) = 1$.
\item Every word in $W_i$ has every word of $L_{i-1,2,0}$ as a subword.
\item Every word $w_{i,j}$ with $0 \leq j \leq n_i-1$ is in $L_{i-1,k,1} \cdot \{\epsilon, 0\}$ for the same $k$, and $w_{i,0} \in L_{i-1,k,1}$.
\end{itemize}

The first two items mean, intuitively, that $W_i$ forms a cycle under Hamming distance, except that going around the cycle once has a cocyclic effect of moving the word to the left by one step. This is where the connectedness of the language will come from. The third item will be important in the proof of minimality. The importance of the fourth condition is that it prevents the accumulation of zeroes when we perform the inductive construction.

Suppose these properties hold for $i$. We build $W_{i+1}$ as follows: First, pick
\[ w_{i+1,0} = w_{i,0} \cdot 0 u \cdot w_{i,0} \]
where $u$ is any sufficiently generic word in $L_{i,k,0} = W_i^k$ for large enough $k$. More precisely, it is sufficient that $u = w_{i,h_0} w_{i,h_1} \cdots w_{i,h_{k-1}}$ satisfies that for all $j, j' \in \{0,\ldots,n_i-1\}$ there exist $|n' - n| \geq 2$ such that $h_n = j, h_{n+1} = j'$ and $h_{n'} = j, h_{n'+1} = j'$.

Consider the decomposition $w_{i+1,0} = w_{i,0} \cdot 0u \cdot w_{i,0}$, and rewrite the central
\[ 0u = 0 w_{i,h_0} w_{i,h_1} \cdots w_{i,h_{k-1}} \]
(where $w_{i,h_j} \in W_i$ for all $j \in \{0, \ldots, k-1\}$) successively into
\[ 0 w_{i,h_0+1} w_{i,h_1} \cdots w_{i,h_{k-1}} \]
\[ 0 w_{i,h_0+2} w_{i,h_1} \cdots w_{i,h_{k-1}} \]
\[ \ldots \]
\[ 0 w_{i,n_i-1} w_{i,h_1} \cdots w_{i,h_{k-1}} \]
\[ w_{i,0} 0 w_{i,h_1} \cdots w_{i,h_{k-1}} \]
\[ w_{i,1} 0 w_{i,h_1} \cdots w_{i,h_{k-1}} \]
\[ \ldots \]
\[ w_{i,h_0} 0 w_{i,h_1} \cdots w_{i,h_{k-1}} \]
and continue by similarly rotating the word $w_{i,h_1}$ `around $W_i$' to move it to the left in $n$ rewriting steps, then $w_{i,h_2}$ and so on until $0u$ has been changed to $u0$. This rewriting is done between the two occurrences of $w_{i,0}$ (the prefix and suffix of the initial word $w_{i+1,0}$). Collecting the words we see during this rewriting process into a list, we have defined $w_{i+1,j}$ for $j \in \{0, 1, \ldots, n_i k\}$.

Next, consider the decomposition $w_{i+1,n_i k} = w_{i,0} \cdot u \cdot 0 w_{i,0}$ and rewrite the prefix $w_{i,0}$ into $w_{i,n_i-1}$ in $n_i-1$ steps. Then rewrite the final word $0 w_{i,0}$ into $w_{i,0} 0$ in $n_i$ steps. Letting $n_{i+1} = n_i (k + 2)$, we have defined $w_{i+1,j}$ for all $j \in \{0,1,\ldots,n_{i+1}-1\}$, and the final word is $w_{i+1,n_{i+1}-1} = w_{i,n_i-1} u w_{i,0} 0$.


By construction (and induction) we have $H(w_{i+1,j}, w_{i+1,j+1}) = 1$ for all $j \in \{0,\ldots,n_{i+1}-2\}$. We have $H(w_{i+1,0} 0, 0 w_{i+1,n_{i+1}-1}) = 1$ because
\[ H(w_{i,0} 0 u w_{i,0} \cdot 0, 0 \cdot w_{i,n_i-1} u w_{i,0} 0) = H(w_{i,0} 0, 0 w_{i,n_i-1}) = 1 \]
by induction.

Every word in $W_{i+1}$ has every word in $L_{i,2,0}$ as a subword because $u$ has at least two copies of each such word, and at any time during the construction of the words $w_{i+1,j}$ we are modifying only one $W_i$-subword of $u$ (and the rest have only been shifted, as they have been fully cycled through). We also have that every word $w_{i+1,j}$ with $0 \leq j < n_{i+1}-1$ is in $L_{i,k,1} \cdot \{\epsilon, 0\}$ directly by construction.

This concludes the construction of the sets $W_i$, and the proof of the inductive properties listed above.

Next, we prove by induction a property which we call property A: for any $i$ and $m \geq i + 1$, we have $W_m \subset L_{i,k,2} \cdot \{\epsilon, 0\}$ for some $k$, and further $w_{m,0} \in L_{i,k,2}$. For $m = i + 1$ this is direct from the last condition in the list of properties already proved, because $L_{i,k,1} \subset L_{i,k,2}$. Now suppose the assumptions hold for $m$. Then when building the words in the next stage, note that we always have a decomposition into a concatenation of words in $W_m$ with either $\epsilon$ or $0$ between them. The words in $W_m$ have a decomposition into $W_i$-words and $\epsilon$, $0$ or $00$ between them by induction. Thus, we immediately see that every word $w \in W_{m+1}$ has a decomposition into $W_i$-words and (a priori) gaps $\epsilon$, $0$, $00$ or $000$ between them. To show property A, we need to consider this decomposition in more detail and show that $000$ does not appear, the gap at the end is of length at most one, and there is no gap at the end of the decomposition if $w = w_{m+1,0}$.

First, any gap in the decomposition of $w$ which occurs inside the decomposition of a $W_m$-word is, by induction, of length at most $2$. Any gap properly inside $w$ and between two $W_m$-words has length at most $2$ since the $W_i$-decomposition of a word in $W_m$ begins with a word of $W_i$ and ends in a gap of length at most~$1$.

The gap at the end of $w$ in its $W_i$-decomposition is precisely as long as the gap in the last $W_m$-word, if $w \neq w_{m,n_m-1}$, since a $0$ appears at the end of the $W_m$-decomposition of $w$ only in this case. Thus in these cases the gap at the end is of length at most one by induction. In $w_{m+1,n_{m+1}-1}$ the gap is one more than the gap after the last $W_i$-word in the decomposition of $w_{m,0}$. By induction, the gap is of length $1$ in this case. We have shown that $w \in L_{i,k,2} \cdot \{\epsilon, 0\}$ for all $w \in W_{m+1}$. Finally, in the decomposition of $w_{m+1,0}$, the last word is $w_{m,0}$, and thus by induction the last word is in $W_i$, so $w_{m+1,0} \in L_{i,k,2}$.

We define our subshift to be $X = \B^{-1}(\bigcup_i W_i)$, where we note that by construction the language $\bigcup_i W_i$ is extendable. Note that for any $i$, the subshift $X$ is contained in the SFT $X_i = \B^{-1}((W_i \{\epsilon, 0, 00\})^*)$ by property A. Note also that it is nonempty since $W_0 \neq \emptyset$ and since $L \subset \B(\B^{-1}(L))$ for an extendable language $L$.

We claim that $X$ is minimal. To see this, let $w$ occur in some point in $X$. Then $w$ is a subword of some $w_{i,j}$ by definition. Observe that $w_{i,j}$ is a subword of every word in $W_{i+1}$. Since $X \subset X_{i+1}$, actually every point in $X$ contains $w_{i,j}$ with bounded gaps, thus $w$ also.

To see that $X$ is language-connected, it is enough to show that for arbitrarily large $n$, the language $\B_n(X)$ is connected, as paths between shorter words are obtained as projections of paths between arbitrary extensions into longer words. Let $n = n_i$ for any $i$ and consider any word $w \in \B_n(X)$.

We first show that there is a path from $w$ to some word in $W_i$. To see this, observe that by the definition of $X$, $w$ appears as a subword of some $W_m$. Thus, it appears as a subword of the $u$-part of $w_{m+1,0} \in W_{m+1}$, at distance at least $\ell_m$ from the boundaries of $u$ (since $u$ contains at least two copies of each word in $L_{m,2,0}$; see the definition of $W_{m+1}$). The path from $w_{m+1,0}$ to $w_{m+1,n_{m+1}-1}$ gives a path from $0u$ to $u0$ by projection.

Consider now the splitting of $w_{m+1,0}$ into a word of $L_{i,k,2} \cdot (\epsilon, 0)$ given by property A. In this splitting, $w$ appears between some two $W_i$-words in the $u$-part of $w_{m+1,0}$. Thus is it a subword of $u$ which is contained in some $v 0^a v'$ for some $v, v' \in W_i$ and $a \in \{0,1,2\}$. Suppose this subword begins in coordinate $c$, that is, $w = (v 0^a v')_{[c,c+n-1]}$. If $w$ is not the suffix of $v 0^a v'$, then the path from $0u$ to $u0$ restricts to a path from $w$ to $(v 0^a v')_{[c+1,c+n]}$, then from $(v 0^a v')_{[c+1,c+n]}$ to $(v 0^a v')_{[c+2,c+n+1]}$, and so on, and finally we connect $w$ to $v' \in W_i$.

The set $W_i$ is connected by construction, so we have shown that the language is connected.
\end{proof}

In the construction, we have more or less complete freedom in the choice of $u$, and it is easy to add some additional nice properties, in particular the subshift can be made uniquely ergodic, to have positive entropy or to have zero entropy. However, we do not know when it is topologically mixing.

\section*{Acknowledgements}

We thank Nishant Chandgotia for discussions and comments.

\bibliographystyle{plain}
\bibliography{bib}{}

\end{document}